\documentclass{amsart}
\usepackage{amsfonts}
\usepackage{amsmath}
\usepackage{amsfonts}
\usepackage{latexsym}
\usepackage{amssymb}
\usepackage[latin1]{inputenc}

\newcommand{\as}{a_1,\ldots,a_n}
\newcommand{\bs}{b_1,\ldots,b_n}

\newcommand{\U}{\mathcal{U}}
\newcommand{\um}{\U_{\cM}}
\newcommand{\vf}{\tau}
\newcommand{\spa}{\sigma({\bar{a}})}
\newcommand{\spb}{\sigma({\bar{b}})}

\newcommand{\ain}{(a_i)_{i=1}^n}

\newcommand{\bin}{(b_i)_{i=1}^n}

\def\refe#1{(\ref{#1})}
\def\ta{\bar{a}}
\def\tb{\bar{b}}
\def\dsm{DS(\cM)}

\newcommand{\RR}{\mathbb{R}}
\newcommand{\NN}{\mathbb{N}}
\newcommand{\MM}{\mathbb{M}}
\newtheorem{fed}{Definition}[section]
\newtheorem{teo}[fed]{Theorem}
\newtheorem{lema}[fed]{Lemma}
\newtheorem{coro}[fed]{Corollary}
\newtheorem{pro}[fed]{Proposition}
\newtheorem{rema}[fed]{Remark}

\newtheorem{exam}[fed]{Example}

\def\mp{M_+ ^{\sim} (\mathbb R^n)}
\def\co{\mathrm{conv}}
\def\re{\mathrm{Re}}
\def\im{\mathrm{Im}}
\def\M{\mathcal{M}}

\def\cA{\mathcal{A}}
\def\cB{\mathcal{B}}

\def\cD{\mathcal{D}}

\def\cM{\mathcal{M}}
\def\cP{\mathcal{P}}
\def\cS{\mathcal{S}}

\def\cU{\mathcal{U}}
\newcommand{\supp}{\mbox{supp }}
\newcommand{\ds}{\displaystyle}
\begin{document}

\title{ The local form of doubly stochastic maps
and joint majorization in II$_1$ factors}
\date{}
\author{Martín Argerami}
\author{Pedro Massey}
\address{Department of Mathematics and Statistics,
    University of Regina, Saskatchewan,Canada
        \ S4S 0A2\\{\tt argerami@math.uregina.ca}}
\address{Departamento de Matem\'atica, Facultad de Ciencias Exactas\\ Universidad Nacional de La Plata,
Argentina\\{\tt massey@mate.unlp.edu.ar}}

\thanks{Supported in part by the Natural Sciences and Engineering Research Council of Canada}
\dedicatory{{\it Dedicated to our families}}
\begin{abstract}
We find a description of the restriction of doubly stochastic maps
to separable abelian $C^*$-subalgebras of a II$_1$ factor $\cM$.
We use this local form of doubly stochastic maps to develop a
notion of joint majorization between $n$-tuples of mutually
commuting self-adjoint operators that extends those of Kamei (for
single self-adjoint operators) and Hiai (for single normal
operators) in the II$_1$ factor case. Several characterizations of
this joint majorization are obtained. As a byproduct we prove that
any separable abelian $C^*$-subalgebra of $\cM$ can be embedded
into a separable abelian $C^*$-subalgebra of $\cM$ with diffuse
spectral measure.
\end{abstract}

\maketitle

\section{Introduction}
Majorization between self-adjoint operators in finite factors was
introduced by Kamei \cite{Ka} as an extension of Ando's
definition of majorization between self-adjoint matrices
\cite{Ando}, a useful tool in matrix theory. Later on, Hiai
 considered majorization in semifinite
factors between self-adjoint and normal operators
\cite{Hiai0,Hiai}. The reason why majorization has attracted the
attention of many researchers (see the discussion in \cite{Hiai}
and the references therein) is that it provides a rather subtle
way to compare operators and occurs naturally in many contexts
(for example \cite{JD,Fack,Doug}). Recently, majorization has
regained interest because of its relation with norm-closed unitary
orbits of self-adjoint operators and conditional expectations onto
abelian subalgebras
\cite{JD,JMD,Arv-Kad,Doug,Kad01,Kad02,neu,she}. One of the goals
of this paper (section \ref{aplic}) is to obtain an extension of
the notion of majorization between normal operators to that of
\emph{joint majorization} between $n$-tuples of commuting
self-adjoint operators in a II$_1$ factor (such extension is
achieved in \cite{Luis} for finite dimensional factors). In order
to obtain characterizations of this extended notion we describe
the \emph{local form} of a doubly stochastic map (DS), i.e. we get
a family of particularly well behaved DS maps that approximate the
restriction of any DS map to separable abelian $C^*$-subalgebras
of a II$_1$ factor (section \ref{sec 3}). As a byproduct, we
construct separable abelian diffuse refinements of separable
abelian $C^*$-subalgebras of a II$_1$ factor $\M$. This
construction seems to have interest on its own. Some of the
techniques we use seem to be new, even in the single element case.

So far we have restricted our attention to the II$_1$ factor case
because, on one hand, technical aspects of the work become simpler
and on the other hand, this is the context where majorization has
its full meaning. 
Since every finite von Neumann algebra acting on a separable
Hilbert space has a direct integral decomposition in terms of
finite factors, the study of II$_1$ factors provides useful
information about more general algebras.


The paper is organized as follows. In section \ref{prel} we recall
some facts about abelian $C^*$-subalgebras of a II$_1$ factor. In
section \ref{sec 3}, after describing some technical results, we
obtain a description of the local structure of doubly stochastic
maps. In section \ref{aplic} we introduce and develop the notion
of joint majorization between finite abelian families of
self-adjoint operators in a II$_1$ factor and we obtain several
characterizations of this relation. Finally, in section
\ref{seccion refinamientos} we prove the results described in
section \ref{sec 3}.

\section{Preliminaries}\label{prel}

Throughout the paper $\cM$ will be a II$_1$ factor with normalized faithful
normal trace $\tau$. The $C^*$-subalgebras of
$\M$ are always assumed unital. The subspace of self-adjoint
elements of $\cM$ will be denoted by $\cM_{sa}$, and we will consider
\textbf{abelian families} $(\as)=(a_i)_{i=1}^n$ in $\cM_{sa}$,
that is finite families of mutually commuting self-adjoint
operators in $\cM$. If $(a_i)_{i=1}^n\subseteq \M_{sa}$ is an
abelian family then $C^*(\as)$ denotes the (unital) separable
abelian $C^*$-subalgebra of $\M$ generated by the elements of the
family. If $\cA$ is an arbitrary abelian $C^*$-subalgebra of $\M$
then $\Gamma(\cA)$ denotes its space of characters, i.e. the set
of *-homomorphisms $\gamma:\cA\rightarrow \mathbb{C}$, endowed
with the weak$^*$-topology. The set $\Gamma(\cA)$ is a compact
space and $\cA\simeq C(\Gamma(\cA))$, where $C(\Gamma(\cA))$
denotes the $C^*$-algebra of continuous functions on
$\Gamma(\cA)$.

\subsection{Joint spectral measures and joint spectral distributions}\label{subseccion del espec conj}

As we will consider a several-variable version of functional calculus,
 we state a few facts about it (see
\cite{tak} for a different description).
 Let $\ta=(a_i)_{i=1}^n$ be an abelian family in
 $\cM_{sa}$. If $\cA=C^*(\as)$, then $\Gamma(\cA)$ can be
embedded in $\prod_{i=1} ^n \sigma(a_{i})\subseteq \mathbb{R}^n$.
In fact, the map $\Phi:\Gamma(\cA)\rightarrow \prod_{i=1} ^n
\sigma(a_{i})\subseteq \mathbb{R}^n$ given by
$\Phi(\gamma)=(\gamma(a_1),\ldots,\gamma(a_n))$ is a
continuous injection and therefore $\Gamma(\cA)$ is homeomorphic
to its image under this map; this image is called the {\bf joint
spectrum} of the family and we denote it by
$\spa\subseteq\prod_{i=1} ^n \sigma(a_{i})$. Note that $\cA\simeq
C(\spa)$ as $C^*$-algebras. If $f\in C(\spa)$, there exists a
normal operator, denoted $f(\as)$, that corresponds to $f$ under the
isomorphism $\cA\simeq C(\spa)$. This association extends the
usual one variable functional calculus.

If $\cA\subseteq \M$ is a separable $C^*$-subalgebra then
$\Gamma(\cA)$ is metrizable and the representation
$C(\Gamma(\cA))\simeq \cA\subseteq \M$ induces a spectral measure
$E_\cA$ \cite[IX.1.14]{con} that takes values on the lattice
$\cP(\cM)$ of projections of $\M$. Let $\mu_\cA$ be the (scalar)
regular Borel measure on $\Gamma(\cA)$ defined by $$
\mu_\cA(\Delta)=\tau(E_\cA(\Delta)).$$ The regularity of $\mu_\cA$
follows from the fact that every open set is $\sigma$-compact
\cite[2.18]{rud}. The map
$\Lambda:L^\infty(\Gamma(\cA),\mu_\cA)\rightarrow\cM$ given by
$\Lambda(h)=\int_{\Gamma(\cA)}\,h\,dE_\cA$ is a normal
$^*$-monomorphism (note that in this case the weak$^*$ topology of
$L^\infty(\Gamma(\cA),\mu_\cA)$, restricted to the unit ball,  is
metrizable) and we have
\begin{equation}\label{integral=traza}
\tau\left(\Lambda(h)\right)=\int_{\Gamma(\cA)}
h\,d\mu_\cA, \ \ \ \forall h\in L^\infty(\Gamma(\cA),\mu_\cA).
\end{equation}
 We consider the von Neumann algebra
$L^\infty(\cA):=\Lambda(L^\infty (\Gamma(\cA),\mu_\cA))
\subseteq \M$.

When $\cA=C^*(a_1,\ldots,a_n)$,
 $E_{\ta}:=E_\cA$ and  $\mu_{\ta}:=\mu_\cA$
are the \textbf{joint spectral measure}
 and  \textbf{joint spectral distribution} of the abelian family $\ta$
and we denote by
$\Lambda_{\ta}:L^\infty(\Gamma(\ta),\mu_{\ta})\rightarrow
L^\infty(\cA)$ the normal isomorphism defined above. It is
straightforward to verify that $\Lambda_{\ta}(\pi_i)=a_i$, $1\leq
i\leq n$, and we write $ h(a_1,\ldots,a_n):=\Lambda_{\ta}(h)$.
 In the case of a single self-adjoint operator $a\in
\M_{sa}$ the measure $\mu_a$ is the usual spectral
distribution of $a$ (see \cite{Arv-Kad}),
and it agrees with the Brown measure of $a$.

In the particular case when $x\in\cM$ is a normal operator, the
real and imaginary parts
 of $x$ are mutually commuting
self-adjoint elements of $\cM$. Identifying the complex plane with
$\RR^2$ in the usual way, it is easy to see that the spectrum of
$x$ as a normal operator coincides with the joint spectrum of the
abelian pair $(\re(x),\im(x))$, and that the spectral measure of
$x$ coincides with the joint spectral measure of
$(\re(x),\im(x))$.

\subsection{Comparison of measures and diffuse measures}

We denote by $\mp$ the set of all regular finite positive Borel
measures $\nu$ on $\mathbb R ^n$ with $\int \|\zeta\|
\,d\nu(\zeta)<\infty$. We write $\nu(f)=\int_{\RR^n} f\ d\nu$, for
every $\nu\in \mp$ and every $\nu$-integrable function
$f$. In what follows, $1$ denotes the
constant function and $\pi_i:\RR^n\rightarrow\RR$ denotes the
projection onto the $i^{th}$ coordinate.

\begin{fed}\label{la defi} We say that $\mu$
is majorized by $\nu$, and we write $\mu\prec \nu$, if for every
$\mu_1,\ldots,\mu_m\in \mp$ with $\sum_{i=1}^m \mu_i=\mu$ there
exist $\nu_1,\ldots,\nu_m\in \mp$ such that $\sum_{i=1}^m \nu_i= \nu$,
$\nu_i(1)=\mu_i(1)$ and $\nu_i(\pi_j)=\mu_i(\pi_j)$ for
$1\leq i\leq m$, $1\leq j\leq n$.
\end{fed}

The relation $\prec$ in Definition \ref{la defi} does not seem to be
called ``majorization'' in the literature, but it will be a
suitable name for us in the light of Theorem \ref{teorema
equivalencias}. If $\mu,\,\nu\in\mp$ we shall write $\nu\sim\mu$
whenever $\nu(1)=\mu(1)$ and $\nu(\pi_j)=\mu(\pi_j)$ for every
$1\leq j\leq n$.

\begin{teo}\label{teo de choquet}
\cite[I.3.2]{alf} Let $\mu,\,\nu\in\mp$. Then  $\mu\prec\nu $
if and only if $\mu(f)\le\nu(f)$ for every continuous convex function
$f:\RR^n\rightarrow \RR$.
\end{teo}

The next corollary is a direct consequence of Theorem \ref{teo de
choquet} and the identity in equation \refe{integral=traza}.
\begin{coro}\label{mayo medidas sii convex functions}
Let $\ta=\ain,\,\bar{b}=\bin\subset\cM_{sa}$ be two abelian families. Then
$\mu_{\ta}\prec\,\mu_{\tb}$ if and only if
$\vf(f(\as))\le\vf(f(\bs))$ for every continuous convex function
$f:\RR^n\rightarrow \RR$.
\end{coro}

We end this section with the following elementary fact about
diffuse (scalar) measures, i.e. measures without atoms (recall
that $x$ is an atom of a measure $\mu$ if $\mu(\{x\})>0$).

\begin{lema}\label{medida difusa}
Let $K\subset\RR^n$ be compact and let $\mu$ be a regular
diffuse Borel probability measure
 on $K$.  Then for every  $\alpha\in(0,1)$ there exists a measurable set
 $S\subset K$ such that $\mu(S)=\alpha$.
\end{lema}

\section{The local form of doubly stochastic maps}\label{sec 3}

A linear map $\Phi:\M\rightarrow \M$ is said to be {\bf doubly
stochastic} \cite{Hiai0} if it is unital, positive, and trace
preserving. We denote the set of all doubly stochastic maps on
$\cM$ by DS$(\cM)$. Doubly stochastic maps play an important role
in the theory of majorization between self-adjoint operators (see
for instance \cite{au1,au2,Hiai0,Hiai}); thus, the study of their
structure appears as a natural topic here.

In what follows we introduce some terminology, we state Theorem
\ref{Adifusion de la medida conjunta2}, Proposition
\ref{Apromedio}, and Lemma \ref{Adixmier} and then we use them to prove
Theorem \ref{la implicacion complicada}. The proofs of these
results will be presented at the end of the paper, in section
\ref{seccion refinamientos}. Although technical, they seem to
have some interest on their own.

Let $\cA\subseteq \M$ be an abelian $C^*$-subalgebra, and let
$E_{\cA}$ and  $\mu_\cA$ denote the spectral measure and the spectral distribution
of $\cA$ as defined in section
\ref{subseccion del espec conj}. If $x\in\Gamma(\cA)$ is such that
$E_\cA(\{x\})\neq 0$, we say that $x$ is an atom for $E_\cA$. The
set of atoms of $E_\cA$ is denoted At$(E_{\cA})$.
 Since
$\mu_{\cA}=\tau\circ E_{\cA}$, the faithfulness of the trace
implies that At$(\mu_{\cA})=$At$(E_{\cA})$. We say that $\cA$ is
\textbf{diffuse} if At$(E_{\cA})=\emptyset$. The following theorem
states that spectral measures of a separable $\cA$ can be refined
in a coherent way.

\begin{teo}\label{Adifusion de la medida conjunta2}
Let $\cA\subseteq \M$ be a separable abelian $C^*$-subalgebra.
Then there exists $a\in \M_{sa}$ such that $C^*(\cA,a)$ is abelian
and  diffuse.
\end{teo}

\bigskip

Since the atoms of $E_\cA$ are in correspondence with
the set of minimal projections of $L^\infty(\cA)$, Theorem
\ref{Adifusion de la medida conjunta2} provides a way to embed
$\cA$ into a separable $C^*$-subalgebra $\tilde \cA=C^*(\cA,a)$
such that $L^\infty(\tilde \cA)$ has no minimal projections (see
Remark \ref{rem dif} for further discussion).

\begin{pro}\label{Apromedio}
Let $\cB\subset \cM$ be a separable, diffuse, and abelian
$C^*$-sub\-algebra.
 Then there
exists an unbounded set $\MM\subseteq \NN$ such that for every
$m\in \MM$ there exist $k=k(m)$ partitions of the unity
$\{q_i^{t,m}\}_{i=1}^m\subseteq \cB'\cap\cM$, $1\leq t\leq k$,
with $\tau(q_i^{t,m})=1/m$ {\rm (}$1\leq i\leq m,\ 1\leq t\leq
k${\rm )}, and such that for each $b\in\cB$, if we let $
\beta_i^{t,m}=m\,\tau(b\,q_i^{t,m})$, then
\begin{equation}\label{aproxi con prom de part}
\lim_{m\rightarrow \infty} \left\|\ b-\frac1k\sum_{t=1}^k
\left(\sum_{i=1}^m\beta_i^{t,m}\,q_i^{t,m}\right)\right\|=0.
\end{equation}
\end{pro}

\bigskip

\begin{rema}\label{reduce}
{\rm For fixed $m$ and partitions of the unity $\{q_i^t\}_{i=1}^m$
$1\leq t\leq k$, the linear map $$\ds b\mapsto
\frac1k\sum_{t=1}^k\left(\sum_{i=1}^m\,m\,\tau(b\,q_i^t)\,q_i^t\right)$$
is a contraction with respect to the operator norm.}
\end{rema}

 We
denote by $\cD(\cM)$ the convex semigroup
$\cD(\cM)=\co\{Ad\,u:\,u\in\cU(\cM)\}$.

\bigskip

\begin{lema}\label{mat doble estoc}
Let $\{p_i\}_{i=1}^m,\,\{q_i\}_{i=1}^m\subseteq\M$ be partitions
of the unity such that $\tau(p_i)=\tau(q_i)=\frac{1}{m},$ and let
$T\in \dsm$. Then there exists $\rho\in \cD(\cM)$ such that if
$\beta_1,\ldots,\beta_m\in\RR$ and $\alpha_i=m\,\sum_{j=1}^m
\beta_j\tau(T(q_j)\,p_i)$ for $1\leq i\leq m$, we have
\begin{equation}\label{estar en la co de la lo} \sum_{i=1}^m
\alpha_i p_i=\rho\left(\sum_{i=1}^m \beta_i \,q_i\right).\end{equation}
\end{lema}
\proof  Let $\gamma_{i,j}=m\,\tau(T(q_j)\,p_i)\ge0$; it is then
straightforward to verify that $(\gamma_{i,j})\in \RR^{m\times m}$
is a doubly stochastic matrix and, moreover, that
$\alpha_i=\sum_{j=1}^m \gamma_{i,j}\,\beta_j$ for every
$i=1,\ldots,m$. By Birkhoff's theorem the doubly stochastic matrix
$(\gamma_{i,j})$ can be written as a convex combination of
permutation matrices, i.e. $(\gamma_{i,j})=\sum_{\sigma\in \mathbb
S_m} \eta_\sigma P_\sigma$, where $\eta_\sigma \geq
0,\,\sum_{\sigma\in \mathbb S_m}\eta_\sigma=1$ and $P_\sigma$ is
the $m\times m$ permutation matrix induced by $\sigma \in \mathbb
S_m$. Then we have
\begin{equation}\label{los alfai}
\alpha_i=\sum_{j=1}^m\gamma_{i,j}\,\beta_j=\sum_{\sigma\in \mathbb
S_m}\eta_\sigma \,\beta_{\sigma(i)}\ \ \ 1\leq i\leq m.
\end{equation}
The fact that $\M$ is a II$_1$ factor and that the elements of the
partitions $\{p_i\}_i$, $\{q_i\}_i$ have the same trace guarantees
the existence of unitaries $u_\sigma$ such that $u_\sigma
\,q_{\sigma(i)}\,(u_\sigma)^*=p_i$, $1\leq i\leq m$, for every
$\sigma\in \mathbb S_m$. Indeed, if $\sigma\in \mathbb S_m$, the
equalities, $\tau(q_{\sigma(i)})=\tau(p_i)$ imply that there exist
partial isometries $v_{i,\sigma}\in \M$ such that
$v_{i,\sigma}v_{i,\sigma}^*=p_i$ and
$v_{i,\sigma}^*v_{i,\sigma}=q_{\sigma(i)}$ for $i=1,\ldots,m$.
Then $u_\sigma=\sum_{i=1}^m v_{i,\sigma}\in \M$ are the required
unitaries.
 Using equation (\ref{los alfai}), and letting $\rho(\cdot\,)=\sum_{\sigma\in\mathbb
S_m}\eta_\sigma \,u_\sigma\,(\cdot\,)\,u_\sigma^*$,
\begin{eqnarray*}\label{del claim}
\sum_{i=1}^m\alpha_i\,p_i&=& \sum_{i=1}^m \left( \sum_{\sigma\in
\mathbb S_m}\eta_\sigma \,\beta_{\sigma(i)}\right) \,p_i=
\sum_{\sigma\in \mathbb S_m}\eta_\sigma  \left( \sum_{i=1}^m
\beta_{\sigma(i)}\, u_\sigma
\,q_{\sigma(i)}\,u_\sigma^*\right)\\ &=& \sum_{\sigma\in \mathbb
S_m}\eta_\sigma  \ u_\sigma\left( \sum_{i=1}^m \beta_{i}\,
 q_{i}\right)u_\sigma^*=\rho\left(\sum_{i=1}^m \beta_{i}\,q_i\right).\qed
\end{eqnarray*}

\bigskip

\begin{lema}\label{Adixmier}
Let $\cB\subset\cM$ be a separable $C^*$-subalgebra, and let
$\{p_i\}_{i=1}^m\subseteq \cB\,'\cap\cM$ be a partition of the
unity. Then there exists a sequence $\{\rho_i\}_{i\in
\NN}\subset\cD(\cM)$ such that for every $b\in\cB$, if we let
$\beta_i(b)=\tau(b\,p_i)/\tau(p_i)$, then
\begin{equation*}
\lim_{j\rightarrow \infty}\left\|\rho_j(b)-\sum_{i=1}^m\beta_i(b)
p_i\right\|=0.
\end{equation*}
\end{lema}

\bigskip

\begin{teo}\label{la implicacion complicada} Let
$\cA,\, \cB\subseteq\M$ be separable abelian $C^*$-sub\-algebras and
let $T\in \dsm$. Let $\cS$ be the operator subsystem of $\cB$
given by $\cS=T^{-1}(\cA)\cap\cB$. Then there exists a sequence
$(\rho_r)_{r\in \NN}\subseteq \cD(\cM)$ such that
$\lim_{r\rightarrow\infty}\left\|T(b)-\rho_r(b)\right\|=0$ for
every $b\in \cS$.
\end{teo}

\proof  First, note that we just have to prove the theorem for
separable diffuse abelian $C^*$-subalgebras of $\M$; indeed,
assume it holds for such algebras and let $\cA,\,\cB\subseteq \M$
be arbitrary separable abelian $C^*$-subalgebras. Then, by Theorem
\ref{Adifusion de la medida conjunta2} there exist separable
diffuse abelian subalgebras $\tilde\cA$ and $\tilde\cB$ of $\cM$
such that $\cA\subseteq \tilde\cA$ and $\cB\subseteq \tilde \cB$.
Thus we get a sequence $\{\rho_r\}_{r\in\NN}\subseteq \cD$ such
that $\lim_{r\rightarrow \infty}\|T(b)-\rho_r(b)\|=0$, for every
$b\in T^{-1}(\cA)\cap\cB\subseteq T^{-1}(\tilde\cA)\cap\tilde
\cB$. So we assume that $\cA$ and $\cB$ are diffuse.

 By Proposition \ref{Apromedio}, there exists an unbounded set
$\MM\subseteq \NN$ and, for each $m\in \MM$,  $k(m)$ partitions of
the unity $\{q_i^{j,m}\}_{i=1}^m \subseteq\cB\,'\cap\cM$ and
$\{p_i^{j,m}\}_{i=1}^m\subseteq \cA\,'\cap\cM$ (in order to
simplify the notation we avoid the supra-index $m$ and write
$q_i^j,\,p_i^j$), $1\leq j\leq k$, such that for every $b\in
T(\cA)^{-1}\cap \cB$ and every $r\in \NN$, there exists
$m_0(r,b)\in \MM$ such that if $m\geq m_0$ we have
\begin{equation}\label{aproxi para b}
\left\|b-\frac1k\sum_{j=1}^k\left( \sum_{i=1}^m \beta_i^j
q_i^j\right)\right\|<\frac1r\end{equation} and
\begin{equation}\label{aproxi para a}
\left\|T(b)-\frac1k\sum_{j=1}^k\left( \sum_{i=1}^m \alpha_i^j
p_i^j\right)\right\|<\frac1r\end{equation} where
$\beta_i^j=m\,\tau(b\,q_i^j)$, $\alpha_i^j=m\,\tau(T(b)\,p_i^j)$,
$\tau(p_i^j)=\tau(q_i^j)=1/m$, (from the construction of such
partitions it is evident
 that we can assume that both have the same unbounded set $\MM$
 and
 the same  $k(m)$ for every $m\in \MM$). Fix $b\in\cB$.
 Since $\|T\|=1$, it
follows from equation (\ref{aproxi para b}) that
\begin{equation}\label{con los por el piso} \left\|T(b)- \frac1k\sum_{j=1}^k\left(
\sum_{i=1}^m \beta_i^j \,T(q_i^j)\right)\right\|\leq
\frac1r.\end{equation} Applying to \refe{con los por el
piso} the fact that the linear map in Remark \ref{reduce} is
linear and contractive (with $\{p_i^j\}_i$ as the partitions of
the unity), we get
\begin{equation}\label{ni te explico...} \left\|
\frac1k\sum_{j=1}^k \sum_{i=1}^m \alpha_i^j \,p_i^j-
\frac{1}{k^2}\sum_{j=1}^k\left(\sum_{t=1}^k \sum_{i=1}^m
\alpha_i^{j,\,t} \,p_i^{\,t}\right)
  \right\|\leq \frac1r,
\end{equation}
where $\alpha_i^{j,\,t} = m\,\sum_{l=1}^m \beta _l^j\,
\tau(T(q_l^j)p_i^{\,t})$, and $\alpha_i^j$ as defined above. By Lemma \ref{mat doble estoc} there exists
$\rho_{j,\,t}^m\in \cD(\cM)$ such that
\begin{equation}\label{esta ec}\sum_{i=1}^m \alpha_i^{j,\,t} p_i^{\,t}
 =\rho_{j,\,t}^m\left(\sum_{l=1}^m \beta_l^j q_l^j\right),
\ \ \ \  1\leq j,t\leq k.
  \end{equation}
 Using (\ref{aproxi para a}),
 (\ref{ni te explico...}), and (\ref{esta ec}) we
get
\begin{equation}\label{es la ultima} \left\|T(b)-
\frac{1}{k^2}\sum_{j=1}^k \sum_{t=1}^k
\rho_{j,\,t}^m\left(\sum_{l=1}^m\beta_l^j\,q_l^j\right)
  \right\|\leq \frac2r,
\end{equation}
By Lemma \ref{Adixmier} there exist sequences $(\tilde \rho_n^j)_{n\in
\NN}\subseteq \cD(\cM)$, $1\leq j\leq k$, independent of $b$, such
that for every $r\in\NN$ there exists $n_0=n_0(r,b)$ such that if
$n\geq n_0$ then
\begin{equation}\label{agus}
\left\|\sum_{l=1}^m \beta_l^j \,q_l^j -\tilde \rho_n^j(b)\right\|\leq
\frac1r, \ \ 1\leq j\leq k.
\end{equation}
 From (\ref{es la ultima}) and (\ref{agus}), together with the fact that each $\rho\in\cD(\cM)$ is contractive
 we get, for  every  $n\geq n_0(r,b)$
\begin{equation}\label{charly}
\left\|T(b)- \frac{1}{k^2}\sum_{j=1}^k \sum_{t=1}^k
\rho_{j,\,t}^m(\tilde \rho_n^j(b))
  \right\|\leq \frac3r,
\end{equation}
Consider a dense countable subset $\{b_1,b_2,\ldots\}$ of $\cB$.
Now define $n(r)$, $m(r)$ as $n(r)=\max\{n_0(r,b_1),\ldots,n_0(r,b_r)\}$ and
$m(r)=\max\{m_0(r,b_1),\ldots,m_0(r,b_r)\}$ and let
$\rho_{r}:=\frac{1}{k^2}\sum_{j=1}^k \sum_{t=1}^k
\rho_{j,\,t}^{m(r)}\circ\tilde \rho_{n(r)}^j\in\cD(\cM)$, where
$k=k(m(r))$.  Then, from the previous calculations, we see that
$\|T(b_j)-\rho_{r}(b_j)\|<3/r$ whenever $1\leq j\leq r$. Let
$b\in\cB$, and $\epsilon>0$. Then there exists $l\in\NN$ such that
$\|b-b_l\|<\epsilon/3$. If $r>\max\{l,9/\epsilon\}$, then
$\|T(b_l)-\rho_{r}(b_l)\|<\epsilon/3$, and so
$\|T(b)-\rho_{r}(b)\|\le \epsilon.$ \qed

\begin{coro}\label{en caso de flias abelianas} Let $T\in DS(\M)$
and let $(a_i)_{i=1}^n,\,(b_i)_{i=1}^n\subseteq \M_{sa}$ be
abelian families such that $T(b_i)=a_i$ for $1\leq i\leq n$. Then
there exists a sequence $(\rho_r)_{r\in \NN}\subseteq \mathcal D$
such that for $1\leq i\leq n$
$\lim_{r\rightarrow\infty}\|a_i-\rho_r(b_i)\|=0$.
\end{coro}
\proof  Consider $\cA=C^*(\as)$ and $\cB=C^*(\bs)$, which are
separable abelian $C^*$-subalgebras of $\M$. Applying Theorem
\ref{la implicacion complicada} to this algebras we
get a sequence $(\rho_r)_{r\in \NN}\subseteq \mathcal D$ such that
$\lim_{r\rightarrow\infty}\|T(b)-\rho_r(b)\|=0$ for every $b\in
T^{-1}(\cA)\cap\cB$. By our choice, $b_i\in
T^{-1}(\cA)\cap\cB$ and so
$\|T(b_i)-\rho_r(b_i)\|=\|a_i-\rho_r(b_i)\|\xrightarrow[]{r} 0$.
\qed

\section{Doubly stochastic kernels and joint majorization}\label{aplic}

We begin by introducing doubly stochastic kernels, which are a
natural generalization of doubly stochastic matrices. We shall use
them to define joint majorization in analogy with \cite{Luis}.

\begin{fed}\label{ds}
Let $(X,\mu_X)$, $(Y,\mu_Y)$ be two probability spaces. A positive
unital linear map $\nu:L^\infty(Y,\mu_Y)\rightarrow
L^\infty(X,\mu_X)$ is said to be a {\bf doubly stochastic kernel}
if $\int_X\,\nu(1_\Delta)\,d\mu_X=\mu_Y(\Delta)$, for every
$\mu_Y$-measurable set $\Delta\subseteq Y$.
\end{fed}

Doubly stochastic kernels between probability
spaces are norm continuous and normal.

\begin{exam}  {\rm Let $X$ and $Y$ be compact
spaces and let $\mu_X$ and $\mu_Y$ be regular Borel probability
measures in $X$ and $Y$ respectively. Consider $D\in
L^1(\mu_X\times \mu_Y)$ and let $\nu(f)(x)=\int_{Y}
D(x,y)\,f(y)\,d\mu_Y(y)$. Then $\nu:L^\infty(X,\mu_X)\rightarrow
L^\infty(Y,\mu_Y)$ is a doubly stochastic kernel if and only if
$D(x,y)\geq 0$ $(\mu_X\times\mu_Y)$-a.e. and $\int_X
D(x,y)\,d\mu_X(x)=1$ $\mu_Y$-a.e, $\int_Y D(x,y)\,d\mu_Y(y)=1$
$\mu_X$-a.e. In particular, if $\mu_X=\mu_Y$ is a measure with
finite support $\{x_i\}_{i=1}^m$ and such that
$\mu_X(\{x_i\})=\frac{1}{m}$ for $1\leq i\leq m$ then $D$ is a
doubly stochastic kernel if and only if the matrix
$(D(x_i,x_j))_{i,\,j}$ is an $m\times m$ doubly stochastic
matrix.}
\end{exam}


\begin{pro}\label{nucleo doble con espectros}
Let $\ta=\ain$, $\tb=\bin\subseteq \cM_{sa}$ be abelian families.
Then the following statements are equivalent:
\begin{enumerate}
\item There exists $T\in \dsm$ such that $T(b_i)=a_i,\,1\leq i\leq n$.
\item There exists a doubly stochastic kernel $\nu:L^\infty(\spb,\mu_{\tb})\rightarrow
L^\infty(\spa,\mu_{\ta})$ such that $\nu(\pi_i)=\pi_i$, $1\leq
i\leq n$.
\end{enumerate}
\end{pro}

\proof   Assume that $T(b_i)=a_i$, $1\leq i\leq n$, with $T\in
DS(\cM)$. Let $\cA=C^*(\as)$, $\cB=C^*(\bs)$. As $\cM$ is a finite
von Neumann algebra, there exists a conditional expectation
$\mathcal P_{\cA}:M\rightarrow L^\infty(\cA)$ that commutes with
$\tau$. Then $\nu=\Lambda_{\ta}^{-1}\circ \mathcal P_{\cA}\circ
T\circ\Lambda_{\tb}$ is the desired doubly stochastic kernel.
Conversely, let us assume the existence of $\nu$ as in {\it 2}.
Let $\mathcal P_{\cB}:\cM\rightarrow L^\infty(\cB)$ be the
conditional expectation onto $L^\infty(\cB)$ that commutes with
$\tau$. Then define
$T=\Lambda_{\ta}\circ\nu\circ\Lambda_{\tb}^{-1}\circ \mathcal
P_{\cB}\in DS(\M)$. Clearly $T(b_i)=a_i$, $1\leq i \leq n$. \qed


\begin{fed}\label{def mayo}
Let $\ta=\ain$, $\tb=\bin$ be two abelian families in $\cM_{sa}$.
We say that $\ta$ is {\bf jointly majorized} by $\tb$ (and we write
 $\ta\prec\tb$) if there exists a
doubly stochastic kernel $\nu:L^\infty(\spb,\mu_{\tb})\rightarrow
L^\infty(\spa,\mu_{\ta})$ such that $\nu(\pi_i)=\pi_i$,
$1\leq i\leq n$.
\end{fed}

\medskip

 If $(x_1,\ldots,x_n)$ is a finite
family in $\cM$, let $\um(x_1,\ldots,x_n)$ denote the
\textbf{joint unitary orbit} of the family with respect to the
unitary group $\um$ of $\cM$, i.e.
\[
\um(x_1,\ldots,x_n)=\{(u^*x_1u,\ldots,u^*x_nu):\,u\in \um\}.
\]
We shall also consider the convex hull of the unitary orbit of a
family $(x_i)_{i=1}^n$, $$\co(\U_\M(x_i)_{i=1}^n)=\{
(\rho(x_i))_{i=1}^n,\ \rho\in \mathcal D\}.$$
 We denote by $\overline{ \co}(\U_\M(x_i)_{i=1}^n)$,
  $\overline{ \co}^w(\U_\M(x_i)_{i=1}^n)$ and $\overline{
\co}^{1}(\U_\M(x_i)_{i=1}^n)$ the respective closures in the
coordinate-wise norm topology, coordinate-wise weak operator
topology, and coordinate-wise $L^1$ topology.

\medskip

\begin{teo}\label{teorema equivalencias} Let $\ta=\ain,\,\tb=\bin$ be abelian families in
$\cM_{sa}$. Then the following statements are equivalent:
\begin{enumerate}
\item\label{mayo} $\ta$ is jointly majorized by
 $\tb$.
\item\label{norma} $\ta\in \overline{\co}(\mathcal U_{\mathcal
M}(\tb))$. \item\label{eleuno} $\ta\in
\overline{\co}^{\,1}(\mathcal U_{\mathcal M}(\tb))$.
\item\label{wot} $\ta\in \overline{\co}^{\,w}(\mathcal
U_{\mathcal M}(\tb))$.
\item\label{mayo medidas} $\mu_{\ta}\prec\mu_{\tb}$.
\item\label{ultimo momento} There exists a completely positive
map $T\in \dsm$ such that $a_i=T(b_i),\,1\leq i\leq n.$
\item\label{mapa} There exists $T\in \dsm$ such
that $a_i=T(b_i),\,1\leq i\leq n.$
\item\label{convex functions} $\tau(f(\as))\leq \tau(f(\bs))$ for
every continuous convex function $f:\RR^n\rightarrow\RR$.
\end{enumerate}
\end{teo}


\begin{rema}
{\rm
Let $x\in \M$ be a normal operator. Recall (see the
last paragraph of section \ref{subseccion del espec conj}) that there is a
natural way to identify the usual spectral measure of $x$ with
that of the abelian pair $(\text{Re}(x), \text{Im}(x))$. If
$T\in DS(\M)$, then since $T$ is positive $T(x)=y$ if and only if
$T(\text{Re}(x))=\text{Re}(y)$ and $T(\text{Im}(x))=\text{Im}(y)$.
 From these facts and Theorem \ref{teorema
equivalencias}, we see that if $x,\,y\in \M$ are normal operators
then $x\prec y$ in the sense of \cite{Hiai} if and only if
$(\text{Re}(x),\text{Im}(x))\prec (\text{Re}(y),\text{Im}(y))$ in
the sense of Definition \ref{def mayo}.

%
}
\end{rema}

\bigskip

Let
$\mathcal P_\mathcal N$ denote the trace preserving
conditional expectation onto the abelian von Neumann subalgebra
$\mathcal N\subseteq \M$. Using Theorem \ref{teorema equivalencias} we can then obtain a
generalization of Theorem 7.2 in \cite{Arv-Kad}.

\begin{coro}Let $\mathcal N\subseteq \M$ be an abelian von Neumann
subalgebra and let $(b_i)_{i=1}^n\subseteq \M_{sa}$ be an abelian
family. Then $(\mathcal P_\mathcal N(b_i))_{i=1}^n\prec
(b_i)_{i=1}^n$.
\end{coro}

\bigskip

In the remainder of the section we prove the implications needed
to prove Theorem \ref{teorema equivalencias}. The single variable
case of the following lemma can be found in \cite{Hiai}.

\bigskip

\begin{lema}\label{con wot implica doble estoc}
Let $\ta=(a_i)_{i=1}^n,\,\tb=(b_i)_{i=1}^n\subseteq \M_{sa}$ be
abelian families. If $\ta\in \overline\co^{\,w}(\mathcal
U_{\mathcal M}(\tb))$ then there exists a completely positive
$T\in\dsm$ such that $a_i=T(b_i)$, $1\leq i\leq n$.
\end{lema}
\proof  Let $\{(b_1^j,\ldots,b_n^j)\}_{j\in J}\subseteq
\co(\um(\bs))$ such that $b_i
^j\xrightarrow[j]{\text{weakly}}a_i$, $1\leq i\leq n$.
 Then there exists a sequence $(\rho_j)_{j\in J}\subseteq
 \mathcal D$
 such that
$(b_1^j,\ldots,b_n^j)$ $=$ $(\rho_j(b_1), $ $\ldots, \rho_j(b_n)),$ for
every $j\in J$. Note that  $\rho_j$ is a completely positive
doubly stochastic map and the net $\{\rho_j\}_{j\in J}$ is norm
bounded. Therefore this net has an accumulation point in the BW
topology \cite{Arv}, i.e. there exists a subnet (which we still
call $\{\rho_j\}_{j\in J}$) and a completely positive map
$T:\cM\rightarrow\cM$ such that $\rho_j(x)\xrightarrow
[j]{\text{weakly}}T(x)$ if $x\in \M$. By normality of the trace,
$T$ is trace preserving, positive and unital. Since
$\rho_j(b_i)=b_i ^j\xrightarrow [j]{\text{weakly}} a_i,$ we have
$T(b_i)=a_i ,\,1\leq i\leq n.$ \qed

\medskip

\begin{lema}\label{doble stoc implica mayo medidas}
Let $\ta=(a_i)_{i=1}^n$, $\tb=(b_i)_{i=1}^n\subseteq \M_{sa}$ be
abelian families. If $\ta\prec\tb$, then $\mu_{\ta}\prec\mu_{\tb}$.
\end{lema}
\proof  By hypothesis, $\ta\prec\tb$; that is, there exists a
doubly stochastic kernel $\nu:L^\infty(\spb,\mu_{\tb})\rightarrow
L^\infty(\spa,\mu_{\ta})$ such that $\nu(\pi_i)=\pi_i$, $1\leq
i\leq n$. Let $\nu_1,\ldots,\nu_m\in\mp$ with $\sum_{j=1}^m
\nu_j= \mu_{\ta}$. Define measures $\nu\,'_j$ by
$\nu_j'(\Delta)=\nu_j(\nu(1_\Delta))$.
By continuity of
$\nu$, $\nu_j'(f)=\nu_j(\nu(f))$ for every $f\in L^\infty(\spb,\mu_{\tb})$. So
$\nu\,'_j(\pi_i)=\nu_j(\nu(\pi_i))=\nu_j(\pi_i)$, $1\leq i\leq n$ and $1\leq j\leq m$, and
similarly $\nu_j(1)=\nu\,'_j(1)$, so that $\nu_j\,\sim\nu\,'_j$, for $1\leq j\leq m$.
Finally,
$\sum_{j=1}^m\nu\,'_j(\Delta)=\sum_{j=1}^m\nu_j(\nu(1_\Delta))
=\mu_{\ta}(\nu(1_\Delta))=\mu_{\tb}(\Delta)$. Therefore
$\sum_{j=1}^m\nu\,'_j= \mu_{\tb}$. We conclude that
$\mu_{\ta}\prec\mu_{\tb}$. \qed

\begin{lema}\label{medidas implica mayo}
Let $\ta=(a_i)_{i=1}^n,\ \tb=(b_i)_{i=1}^n\subset\cM_{sa}$ be
abelian families. If $\mu_{\ta}\prec\mu_{\tb}$, then there exists
$T\in DS(\M)$ such that $T(b_i)=a_i$, $1\leq i\leq n$. 
\end{lema}
\proof  By compactness, we can consider partitions
$\{\Delta_j^r\}_{j=1}^{m(r)}$ of $\spa$ with
$\mbox{diam}(\Delta_j^r)<1/r$ for every $1\leq j\leq m$. Fix
points $x_1^r,\ldots,x_{m(r)}^r$ with $x_j^r \in\Delta_j^r$ and
define measures $\mu_j^r$ by
$\mu_j^r(\cdot\,)=\mu_{\ta}(\cdot\,\cap\Delta_j^r)$. Then clearly
$\sum_j\mu_j^r=\mu_{\ta}$. As $\mu_{\ta}\prec\mu_{\tb}$ by
hypothesis, there exist measures $\nu_j^r$ with
$\nu_j^r\sim\mu_j^r$ and $\sum_j\nu_j^r=\mu_{\tb}$. Let $g_j^r$ be
the Radon-Nikodym derivatives $g_j^r=d\nu_j^r/d\mu_{\tb}$. Note
that $\sum_jg_j^r=1$ ($\mu_{\tb}-\text{a.e.}$). Define a function
$D_r:\spa\times\spb\rightarrow\RR$ by
\[D_r(s,t)=\sum_{j=1}^{m(r)}\,\frac{g_j^r(t)}
{\mu_{a}(\Delta_j^r)}\,1_{\Delta_j^r}(s).\] We will use the
kernels $D_r$ to approximate $T$. Let us define
$\nu_r:L^\infty(\spb,\mu_{\tb})\rightarrow
L^\infty(\spa,\mu_{\ta})$ by
$$\nu_r(b)(s)=\int_{\spb}\,b(t)\,D_r(s,t)\,d\mu_{\tb}(t).$$ The map
$\nu_r$ can be seen to be doubly stochastic using the equivalence
$\mu_j^r\sim\nu_j^r$. By Proposition \ref{nucleo doble con
espectros} there is an associated sequence $\{T_r\}_r\subset \dsm$
such that $T_r(b_i)=\int _{\spa} \nu_r(\pi_i)\,dE_{\ta}\in
L^\infty(\cA)$, $1\leq i\leq n.$ The bounded net
$\{T_r\}_{r\in\NN}$ has a subnet $\{T_k\}_{k\in K}$ that converges
to a cluster point $T\in \dsm$ in the BW topology. Since this
subnet is bounded, $T(b_i)=w\mbox{-}\lim_{k\in K} T_k (b_i)\in
L^\infty(\cA)$. We claim that $T(b_i)=a_i,\,1\leq i\leq n.$ To see
this, since the net $\{T_k(b_i)\}_{k\in K}$ is bounded, we just
have to prove that $$ \lim_{k} \tau(x\,T_k(b_i))=\tau(x\, a_i),\ \
1\leq i\leq n,\ \ \forall x\in\cA.$$ Equivalently, we have to show
that for every continuous function $f\in C(\spa)$ and every $i=1,\ldots,n$,
 \begin{equation*}
\lim_{k} \int_{\spa} f(s) \left(\int_{\spb}
D_k(s,t)\,\pi_i(t)\,d\mu_{\tb}(t)\right)\,d\mu_{\ta}(s)=\int_{\spa}
f(s)\,\pi_i(s)\,d\mu_{\ta}(s).\end{equation*} This can be seen by
a standard approximation argument, using the uniform continuity of
$f$, the fact that the diameters of $\Delta_j^r$ tend to 0 as $r$
increases, and the equivalence $\mu_j^r\sim\nu_j^r$.
\qed

\bigskip

\noindent{\bf Proof of Theorem \ref{teorema equivalencias}.
}Proposition \ref{nucleo doble con espectros} shows the
equivalence \refe{mapa}$\Leftrightarrow$\refe{mayo} and Corollary
\ref{en caso de flias abelianas} is
\refe{mapa}$\Rightarrow$\refe{norma}.
 The
implication \refe{norma}$\Rightarrow$
\refe{eleuno}$\Rightarrow$\refe{wot} is trivial. Lemma \ref{con
wot implica doble estoc} shows that
\refe{wot}$\Rightarrow$\refe{ultimo momento}, and it is clear that
\refe{ultimo momento}$\Rightarrow$\refe{mapa}.  Lemmas \ref{doble
stoc implica mayo medidas}, \ref{medidas implica mayo} and
Proposition \ref{nucleo doble con espectros} prove the equivalence
 \refe{mayo medidas}$\Leftrightarrow$\refe{mayo}.
So we have that \refe{mayo}-\refe{mapa} are equivalent. Finally,
Corollary \ref{mayo medidas sii convex functions} shows that
\refe{mayo medidas}$\Leftrightarrow$\refe{convex functions}. \qed

\section{Joint unitary orbits of abelian families in $\M_{sa}$}

Given families $\ta=(a_i)_{i=1}^n,\,\tb=(b_i)_{i=1}^n\subseteq \M$,
we say that $\ta$ and $\tb$ are \textbf{jointly approximately
 unitarily equivalent  in $\M$} if $\ta\in \overline{\U_\M(\tb)}$,
that is if there exists a sequence of unitary
operators $(u_n)_{n\in \NN}\subseteq \M$ such that
$\lim_{n\rightarrow \infty}\|u_n b_iu_n^*-a_i\|=0$ for every
$i=1\ldots,n$.  It is clear that this is an equivalence
 relation. Moreover, if $\ta$ and $\tb$ are jointly approximately
 unitarily equivalent in $\M$ then $\ta$ is an abelian family if
 and only if $\tb$ is.
  In \cite{Arv-Kad} a characterization of this equivalence relation
between selfadjoint operators is obtained, in terms of the
spectral distributions. The following results exhibits a list of
 characterizations of this relation for abelian families in
 $\M_{sa}$.

\bigskip

\begin{teo}\label{doble mayo}
Let $\ta=\ain$ and $\tb=\bin\subset\cM_{sa}$ be abelian
families. Then the following statements are equivalent:
\begin{enumerate}
\item\label{1} $\ta$ and $\tb$ are jointly approximately unitary equivalent
in $\M$.
\item \label{2} $\ta\prec \tb$ and $\tb\prec \ta$
\item \label{3} $\mu_{\ta}=\mu_{\tb}$
\item \label{4} $\tau(f(\as))=\tau(f(\bs))$ for every continuous convex function
$f:\RR^n\rightarrow \RR$.
\item \label{5} $\tau(f(\as))=\tau(f(\bs))$ for every continuous function
$f:\RR^n\rightarrow \RR$.
\end{enumerate}
\end{teo}
\proof  By Theorem \ref{teorema equivalencias} we have
(\ref{1})$\Rightarrow$(\ref{2}) and
(\ref{2})$\Leftrightarrow$(\ref{4}). Moreover, (\ref{4}) is
equivalent to $\mu_{\ta}(f)=\mu_{\tb}(f)$ for every convex
function $f$. Then $\mu_{\ta}(f)=\mu_{\tb}(f)$ for every
continuous function $f$ \cite[Proposition I.1.1]{alf}, and this in
turn implies that $\mu_{\ta}=\mu_{\tb}$. Therefore,
(\ref{4})$\Rightarrow$(\ref{5})$\Rightarrow$(\ref{3}). Again, by
Theorem \ref{teorema equivalencias}
(\ref{3})$\Rightarrow$(\ref{2}) and so (\ref{2})-(\ref{5}) are
equivalent. Finally, we prove that
(\ref{3})$\Rightarrow$(\ref{1}). If we assume that
$\mu_{\ta}=\mu_{\tb}$ then,
$\spa=\supp\,\mu_{\ta}=\supp\mu_{\tb}=\spb$ and for every Borel
set $\Delta$ in $\spa$ we have
\begin{equation}\label{igualtraza}
\tau(E_{\ta}(\Delta))=\mu_{\ta}(1_\Delta)=\mu_{\tb}(1_\Delta)=
\tau(E_{\tb}(\Delta)).\end{equation}
Let $\epsilon>0$. By compactness, choose $B_1,\ldots,B_m$ to be a
finite disjoint covering of $\spa=\spb$ such that there are points $x_j\in B_j$
with the property that $|\pi_i(\lambda)-\pi_i(x_j)|<\epsilon/2$
for every $\lambda\in B_j$, $1\leq i\leq n$, $1\leq j\leq m$. Then
we get, using the Spectral Theorem,
\[
\left\|\,a_i-\ds\sum_{j=1}^m\pi_i(x_j)E_{\ta}(B_j)
\,\right\|<\frac{\epsilon}{2}\mbox{, }
\left\|\,b_i-\ds\sum_{j=1}^m\pi_i(x_j)E_{\tb}(B_j)
\,\right\|<\frac{\epsilon}{2}
\] for $i=1,\ldots,n$.
From equation \refe{igualtraza} we get that
$\tau(E_{\ta}(B_j))=\tau(E_{\tb}(B_j))$ for every $j=1,\ldots,m$.
As in the proof of Lemma \ref{mat doble estoc}, we get a unitary $w_\epsilon\in \U(\M)$
such that $w_\epsilon^*E_{\tb}(B_j)w_\epsilon=E_{\ta}(B_j)$ for every
$j$. Then
\[w_\epsilon^*\left(\sum_{j=1}^m\pi_i(x_j)E_{\tb}(B_j)\right)w
_\epsilon=\sum_{j=1}^m\pi_i(x_j)E_{\ta}(B_j).\] Finally, for every
$i$ we have
\[ \|\,w_\epsilon^*b_iw_\epsilon-a_i\|\le\left\|\,w_\epsilon^*\left(b_i
-\ds\sum_{j=1}^m\pi_i(x_j)E_{\tb}(B_j)\right) w_\epsilon\,\right\|
+
\frac\epsilon2<\epsilon.\qed\]

\bigskip

\begin{coro}\label{los iso de medida}
Let $\Theta$ be an automorphism of $\cM$. Then $\Theta|_{\cA}$ is
approximately inner for each separable abelian $C^*$ subalgebra
$\cA\subset\cM$.
\end{coro}

\proof  The uniqueness of the trace guarantees that $\Theta$ is
trace-preserving. Being multiplicative, the range of an abelian
set will be again abelian. So $\Theta$ is a DS map that takes an
abelian family in $\cM$ into another. Consider a countable dense
subset $\{a_i\}$ of $\cA$, and use Theorem \ref{doble mayo} to
obtain unitaries $u_n$ for each
 finite subset $\{a_1,\ldots,a_n\}$. An $\epsilon/3$
argument shows then that the sequence $\{\mbox{Ad}\,u_n\}$
approximates $\Theta$ in all of $\cA$. \qed

\bigskip

Given $\bar{x}=(x_i)_{i=1}^n\subseteq \M$ we denote by
$\overline{\U_\M(\bar x)}^s$ the closure in the coordinate-wise
strong operator topology. An immediate consequence of Theorem
\ref{doble mayo} is that the norm closure of the unitary orbit of
a selfadjoint abelian family in a II$_1$ factor is strongly closed. This
generalizes \cite[Theorem 5.4]{Arv-Kad} and \cite[Theorem
8.12(1)]{she}:

\begin{coro}
Let $\ta=(a_i)_{i=1}^n\subseteq \M_{sa}$ be an abelian family.
Then
$\overline{\U_\M(\ta)}^{\|\,\|}=\overline{\U_\M(\ta)}^s$.
\end{coro}

\proof Let $\tb=(b_i)_{i=1}^n\in\overline{\U_\M(\ta)}^s$.
There exists a net $(b_1^{j},\ldots,b_n^{j})_{j\in J} \subseteq
\cU_\M(\ta)$ such that $b_i^{j}$ converges strongly to $b_i$ for
each $i=1,\ldots,n$. Let $f:\RR^n\rightarrow\RR$ be a continuous
function. Then
$\tau(f(b_1^{j},\ldots,b_n^{j}))=\tau(f(a_1,\ldots,a_n))$ for
every $j$. Using \cite[Lemma II.4.3]{tak} we conclude that
$\tau(f(b_1,\ldots,b_n))=\tau(f(a_1,\ldots,a_n))$. So $(5)$ of
Theorem \ref{doble mayo} implies that
$\tb\in\overline{\cU_\M(\ta)}$. The other inclusion is
trivial.\qed


\section{Some technical results}\label{seccion
refinamientos}

In this section we prove the results presented at the
beginning of section \ref{sec 3}. First, we show that any
separable abelian $C^*$-subalgebra of $\M$ can be embedded into a
 separable diffuse abelian $C^*$-subalgebra. Then, we prove some approximation results that hold
for separable diffuse abelian $C^*$ subalgebras of $\M$.

\subsection{Refinements of spectral measures}

We begin by recalling some elementary facts about inclusions of
abelian $C^*$ algebras. If $\cA\subseteq \cB$ are unital
$C^*$-algebras, then  the function $\Phi:\Gamma(\cB)\rightarrow
\Gamma(\cA)$ given by $\Phi(\gamma)=\gamma|_{\cA}$ is a continuous
surjection onto $\Gamma(\cA)$. If we assume further that
$\cA\subseteq \cB\subseteq \M$ are separable and that
$E_\cA,\,E_\cB$ denote their spectral measures, then
$E_\cA=E_\cB\circ\Phi^{-1}$ and
$\mu_{\cA}=\mu_{\cB}\circ\Phi^{-1}$.

Note that At$(\mu_{\cA})=$At$(E_{\cA})$ where At$(E_{\cA})$ is the
set of atoms of the spectral measure $E_{\cA}$ (see the beginning
of section \ref{sec 3}). Let $\sum_{x\in
\text{At}(E_\cA)}\mu_\cA(\{x\})$ be the \textbf{total atomic mass}
of $E_\cA$. Since $\mu_{\cA}$ is finite, the total atomic mass is
bounded and thus, the set of atoms is countable set.

\begin{lema}\label{la masa atomica decrece}
With the notation above, if $x\in \mbox{At}(E_\cB)$ then
$\Phi(x)\in \mbox{At}(E_\cA)$, and the total atomic mass of $\cB$
is smaller that the total atomic mass of $\cA$.
\end{lema}
\proof Let $x\in\text{At}(E_\cB)$ and note that
$0\neq\mu_\cB(\{x\})\leq
\mu_\cB(\Phi^{-1}(\Phi(\{x\})))=\mu_\cA(\Phi(\{x\})),$ so
$\Phi(x)\in \text{At}(E_\cA)=\text{At}(\mu_\cA)$. We consider the
equivalence relation in $\text{At}(E_\cB)$ induced by $\Phi$, i.e.
$x\sim y$ if $\Phi(x)=\Phi(y)$. If $Q\in \mathcal
Q:=\text{At}(E_\cB)/\sim$ is such that $\Phi(x)=x_Q$ for every
$x\in Q$, then using that $Q$ is countable we get $\sum_{x\in
Q}\mu_\cB(\{x\})=\mu_\cB(Q)\leq
\mu_\cB(\Phi^{-1}(\{x_Q\}))=\mu_\cA(\{x_Q\}).$ Therefore
$$\sum_{x\in \text{At}(E_\cB)} \mu_\cB(\{x\})=\sum_{Q \in \mathcal
Q}\ \sum_{x\in Q}\mu_\cB(\{x\})\leq \sum_{Q\in \mathcal Q
}\mu_\cA(x_Q)\leq \sum_{x\in \text{At}(E_\cA)} \mu_\cA(\{x\}).$$
\qed

\begin{pro}\label{refinamiento de una med spec}
With the notations above, let $x_0\in \Gamma(\cA)$ be an atom of
$E_\cA$ and let $\alpha,\beta\in\RR$ with $0<\alpha<\beta$. Then
there exists $a\in \cA'\cap\M_{sa}$ with $[\alpha,\beta]\subseteq
\sigma(a)\subseteq [\alpha,\beta]\cup\{0\}$,
$P_{\overline{R(a)}}=E_\cA(\{x_0\})$, and such that $E_\cB$ has no
atoms in the fibre $\Phi^{-1}(x_0)$, where $\cB=C^*(\cA,a)\subset\M$.
\end{pro}
\proof  Let $p=E_\cA(\{x_0\})$ and consider a masa $\cA\subset\M$
such that $p\in\cA$. Then $p\cA$ is a masa in the II$_1$ factor
$p\M p$, where the trace is $\tau_p=\frac1{\tau(p)}\,\tau$. It is
well known that there exists a countably generated, non-atomic von
Neumann subalgebra $\cB$ of $p\cA$ such that there is a von
Neumann algebra isomorphism
$\Phi:L^\infty([0,1],m)\rightarrow\cB$, with $m$ the Lebesgue
measure on $[0,1]$, and with $\tau_p(\Phi(f))=\int_0^1\,f\,dm$. Put
$\tilde{a}=\Phi(id)$; it is clear that $\tilde{a}$ has no atoms in its spectrum with the exception of 0, and
that $E_{\tilde a}(\{0\})=1-p$, $\sigma(a)=[0,1]$. Let
$a=(\beta-\alpha)\tilde a+\alpha\,p$, so $[\alpha,\beta]\subseteq
\sigma(a)\subseteq [\alpha,\beta]\cup\{0\}$,
$P_{\overline{R(a)}}=p=E_\cA(\{x_0\})$.  As $p$ is
a minimal projection in $L^\infty(\cA)$, we have
$pb=pbp=\lambda_bp$ for every $b\in \cA$ and so
$ab=apb=\lambda_bpa=bpa=ba$. Thus $a\in\cA'\cap\cM$.

Let $\cB=C^*(\cA,a)$ and let $\Phi:\Gamma(\cB)\rightarrow
\Gamma(\cA)$, $\Psi:\Gamma(\cB)\rightarrow \Gamma(C^*(a))$ be the
continuous surjections induced by the inclusions $\cA\subseteq
\cB$ and $C^*(a)\subseteq \cB$. Note that the restriction
$\Psi|_{\Phi^{-1}(x_0)}$ is injective. Indeed, let
$x,\,y\in\Phi^{-1}(x_0)$ be such that $\Psi(x)=\Psi(y)$, i.e. the
restriction of the characters to $C^*(a)$ coincide. Since
$\Phi(x)=\Phi(y)\,(=x_0)$, the characters also coincide on $\cA$
and therefore are equal as characters in $\cB$, since $\cB$ is
generated by $\cA$ and $C^*(a)$.

On the other hand, if $x\in \Gamma(\cB)$ is such that $x(a)\neq
0$, then $\Phi(x)=x_0$. Indeed, assume that $\Phi(x)\neq x_0$. Let
$f\in C(\Gamma(\cA))$ with $f(\Phi(x))=0$ and $f(x_0)=1$. So
$f\circ\Phi\geq 1_{\Phi^{-1}(x_0)}$. But then $$
\int_{\Gamma(\cB)} f\circ \Phi\
dE_\cB\geq\int_{\Gamma(\cB)}1_{\Phi^{-1}(x_0)}\,d E_{\cB}=
E_\cB(\Phi^{-1}(x_0)=E_\cA(\{x_0\})=p.$$ Note that if
$0\in\sigma(a)$ then it is an isolated point, so in any case we
have $p\in C^*(a)\subseteq\cB$. Then
 $0=f\circ \Phi(x)\geq x(p)\geq 0$,
so $x(p)=0$. Since $0\leq a\leq \beta\,p$, $x(a)=0$ and the claim
follows.

\smallskip

Now let $z\in\Phi^{-1}(x_0)$. If $z(a)\neq 0$, from the first part
of the proof we deduce that $\Psi^{-1}(\Psi(z))=\{z\}$. Therefore
$E_\cB(\{z\})=E_a(\{\Psi(z\}))=0$, since $\Psi(z)(a)\neq 0$ and
At$(E_a)\subseteq \{0\}$. If $z(a)=0$, then
$$\begin{array}{rcl}\{z\}&=&\Phi^{-1}(x_0)\setminus\{x\in\Phi^{-1}(x_0):\
x(a)\neq 0\}\\ &=&\Phi^{-1}(x_0)\setminus\Psi^{-1}(\{x\in
\Gamma(C^*(a)):\ x(a)\neq 0\})\end{array}$$ and
\begin{eqnarray*}
E_\cB(\Psi^{-1}(\{x\in \Gamma(C^*(a)):\ x(a)\neq 0\}))&=&
E_a(\{x\in \Gamma(C^*(a)):\ x(a)\neq 0\})\\
&=&E_\cA(\{x_0\})=E_\cB(\Phi^{-1}(x_0)).
\end{eqnarray*}From this we conclude that $E_\cB(\{z\})=0$.
\qed

\bigskip

\noindent {\it Proof of Theorem \ref{Adifusion de la medida
conjunta2}}. Recall that the set At$(E_\cA)$ of atoms of $E_\cA$
is a (possibly infinite) countable set. If At$(E_\cA)=\emptyset$
then $E_\cA$ is already diffuse and we are done. Otherwise, let us
enumerate At$(E_\cA)=\{x_i:\,1\leq i\leq r\}$, where
$r\in\NN\cup\{\infty\}$. For $1\leq i\leq r$, let
$I_i=[1+\frac{1}{2n},1+\frac{1}{2n-1}]$. Then $I_i\cap
\overline{\bigcup_{1\leq i\neq j\leq r} I_j}=\emptyset$
 and $\bigcup_{i=1}^rI_i\subseteq [1,2]$.
 For each $i=1,\ldots,r$ there exists, by Proposition
\ref{refinamiento de una med spec}, $a_i\in\cA'\cap\M_{sa}$ such
that $P_{\overline{R(a_i)}}=E_{\cA}(\{x_i\})$,
$I_i\subseteq\sigma(a_i)\subseteq I_i\cup\{0\}$,
 and such that $E_{\cA_i}$ has no atoms in the fibre $\Phi_i^{-1}(x_i)$, where
 $\Phi_i:\Gamma(\cA_i)\rightarrow \cA$ denotes the
 continuous surjection induced by the inclusion $\cA\subseteq
 \cA_i:=C^*(\cA,a_i)$.
 Let $a=\sum_{i=1}^r a_i\in\cA'\cap\M_{sa}$
 (this sum converges because the
ranges of the operators $a_i$ are orthogonal and $\|a_i\|\le2$ for
every $i$). Then $\cB:=C^*(\cA,a)$ is an abelian subalgebra of
$\cM$.

 We claim that the spectral measure $E_\cB$ of $\cB$ has no atoms.
 Indeed, first
note that $1_{I_i}\in C(\cup_{1\leq j\leq r}I_j)$ is a continuous
function (because the distance between the sets $I_i$ and
$\cup_{i\neq j}I_j$ is positive); then, since $1_{I_i}(a)=a_i$, it
follows that $\cA_i\subset\cB$ for every $i=1,\ldots,r$.
 Assume now that $x\in
\text{At}(\Gamma(\cB))$ and let $\Phi:\Gamma(\cB)\rightarrow
\Gamma(\cA)$ be as before. By Lemma \ref{la masa atomica
decrece} there exists $i\in\{1,\ldots,r\}$ such that
$\Phi(x)=x_i\in\text{At}(E_\cA)$ . Since $\Phi=\Phi_i\circ\Psi_i$,
where $\Psi_i:\Gamma(\cB)\rightarrow \Gamma(\cA_i)$ is the
surjection induced by the inclusion $\cA_i\subseteq \cB$, we
conclude that $\Psi_i(x)\in\Phi_i^{-1}(x_i)$ is an atom of the
measure $E_{\cA_i}$, again by Lemma \ref{la masa atomica decrece}.
But this last assertion is a contradiction because by construction
there are no atoms in the fibre $\Phi_i^{-1}(x_i)$ by
construction. \qed

\begin{rema}\label{rem dif} {\rm Given an
abelian $C^*$ subalgebra $\cA\subset\cM$, a direct way to find an
abelian $C^*$-subalgebra $\cA\subseteq\tilde\cA\subset\cM$ with
diffuse spectral measure is to consider  a masa in $\cM$ that
contains $\cA$. The additional information we obtain from Theorem
\ref{Adifusion de la medida conjunta2} is that $\tilde \cA$ can be
chosen separable (as a $C^*$-algebra) whenever $\cA$ is separable.
When this is the case, the character space of $\tilde\cA$ is
metrizable, a fact that is crucial for our calculations.}
\end{rema}

\subsection{Discrete approximations in separable diffuse abelian
algebras}

Given a compact metric space it is always possible to find, using
uniform continuity, discrete uniform approximations of continuous
functions by linear combinations of characteristic functions of
certain sets $\{Q_i\}_{i=1}^m$. But if we consider a measure on
this space and we require equal measures for these sets, there
might not be any good uniform approximation based on
characteristic functions (even for measures of compact support in
the real line). Proposition \ref{Apromedio} is an intermediate solution to this problem. It was inspired
by the proof of \cite[Lemma 4.1]{Hiai}.

\bigskip

\noindent {\it Proof of Proposition \ref{Apromedio}.} The space
$\Gamma(\cB)$ is a metrizable compact topological space, so we
consider a metric $d$ in $\Gamma(\cB)$ inducing its topology. Let
$r\in \NN$; by compactness, there exists a partition $\{\tilde
Q_i\}_{i=1}^{k_0}$ of $\Gamma(\cB)$ with diam$_d(\tilde Q_i)<\frac
1 r$ and $\sum_{i=1}^{k_0}\mu_\cB(\tilde Q_i)=1$. Let $m=m(r)$
be such that $1/m\leq \min\{\mu_\cB(\tilde Q_j)^2:\,1\leq j\leq
k_0\}$. Then for $1\leq j \leq
k_0$ there exists $k_j\in\NN$ such that $\mu_\cB(\tilde
Q_j)=k_j/m+\delta_j$ with $0\leq \delta_j< 1/m$. If we let
$\tilde{k}=\tilde{k}(r)=\min_j\{k_j\}$ then $\tilde{k}\geq
\max\{\mu_\cB(\tilde Q_j)^{-1}-1,\,1\leq j\leq k_0\}$.

For $t=1,\ldots,k_0$, choose $\tilde{k}$ partitions
$\{\tilde{Q}_{j,s}^t\}_{s=0}^{k_j}$ of each $\tilde{Q}_j$ ($1\leq
t \leq \tilde{k}$), with $\mu_{\cB}(\tilde{Q}_{j,s}^t) =1/m$ if
$1\leq s\leq k_j$ and $\mu_{\cB}(\tilde{Q}_{j,0}^t)=\delta_j$, in
such a way that $\tilde{Q}_{j,0}^t\subset\tilde{Q}_{j,\,t}^{1}$,
$2\leq t \leq\tilde{k}$. Note that we can always make such a choice: using
Lemma \ref{medida difusa} choose $\tilde{Q}^t_{j,0}\subseteq
\tilde{Q}^1_{j,\,t}$ with $\mu_{\cB}(\tilde{Q}^t_{j,0})=
\delta_j<1/m$, and then take a partition $\{\tilde
Q_{j,s}^t\}_{s=1}^{k_j}$ of $\tilde Q_j\setminus\tilde{Q}^t_{j,0}$
using again Lemma \ref{medida difusa} (note that $\mu_\cB(\tilde
Q_j\setminus\tilde{Q}^t_{j,0})=k_j/m$). By this choice,
$\tilde{Q}^t_{j,0}\cap\tilde{Q}^{t\,'}_{j,0}=\emptyset$ if $t\neq
t\,'$. 

For each $t=1,\ldots,\tilde{k}$, let
$\tilde{Q}_{0,0}^t=\cup_{j=1}^{k_0}\,\tilde{Q}_{j,0}^t$. Then
$\mu_{\cB}(\tilde{Q}_{0,0}^t)=1-\sum_j k_j/m=(m-\sum_{j=1}^{k_0}
k_j)/m$. Finally, make partitions of each set $\tilde{Q}_{0,0}^t$
into $n_1=m-\sum_j k_j$ subsets $\{\tilde Q_i^t\}_{i=1}^{n_1}$ of
measure $1/m$. By re-labeling  the $\tilde{k}$ partitions
$\{\tilde Q_{j,s}^t\}_{j,\,s}\cup\{\tilde Q_i^t\}_i$, we end up
with $\tilde{k}$ partitions $\{Q_i^{t,\,m}\}_{i=1}^m$, for $1\leq
t\leq \tilde{k}$, such that:
\begin{enumerate}
\item[1.] $\mu_\cB(Q_i^{t,\,m})=1/m$, for every $i\in\{1,\ldots,m\}$,
$t\in\{1,\ldots,\tilde{k}\}$;
\item[2.]  diam$_d(Q_i^{t,\,m})\leq 1/r$, if $i>n_1$;
\item[3.] if $1\leq i,\,i\,'\leq n_1$ then $Q_i^{t,\,m}\cap
Q_{i\,'}^{t\,',\,m}=\emptyset$ if $i\neq i\,'$ or $t\neq t\,'$.
\end{enumerate}

Note that the construction of the $k$ partitions
$\{Q_i^{t,\,m}\}_{i=1}^m$ was done in such a way that the subsets
that do not have small diameters are disjoint, even for different
partitions.

Let $\MM=\{m(r),\,r\geq 1\}$ and for every $m=m(r)\in\MM$ let
$k(m)=\tilde{k}(r)$ as defined above and, for $i,t,m$, let
$q_i^{t,m}=E_\cB(Q_i^{t,\,m})$. The set $\MM$ is unbounded
because the measure $\mu_\cB$ being diffuse makes
$\lim_{r\rightarrow \infty}m(r)=\infty$, and so
$\lim_{r\rightarrow \infty}\tilde{k}(r)=\infty$.
For each $t=1,\ldots,k$, $\{q_i^{t,m}\}_{i=1}^m\subset
\cB'\cap\cM$ is a partition of the unity.

Let $b\in \cB$, $\epsilon>0$, and let $f\in C(\Gamma(\cB))$ be
such that $b=\int_{\Gamma(\cB)}f\,dE_\cB$. Then, by compactness,
there exists $\delta>0$ such that if $Q\subseteq \Gamma(\cB)$ with
diam$_d(Q)<\delta$ then diam$(f(Q))<\epsilon$. Let $r\in \NN$ be
such that $1/r<\delta$ and $2\|b\|/k(r)\leq \epsilon$;
 let $m=m(r)\in \MM$, and let
$\beta_i^{t,m}=m\,\tau(b\,q_i^{t,m})=m\int_{Q_i^{t,m}}f\,d\mu_{\cB}$.
Properties 1-3 translate then into
\begin{enumerate}
\item[1'.] $\tau(q_i^{t,m})=1/m$, for every
$i\in\{1,\ldots,m\}$, $t\in\{1,\ldots,k\}$;
\item[2'.] if $i>n_1$, then $|f(x)- \beta_i^{t,m}|\leq \epsilon,\,\forall x\in Q_i^{t,m}$;
\item[3'.] if $1\leq i,\,i\,'\leq n_1$ then $q_i^{t,m}\perp
q_{i\,'}^{t\,',m}$ if $i\neq i\,'$ or $t\neq t\,'$.
\end{enumerate}
Therefore we have
\begin{eqnarray*}
\left\|\,b-\frac1k\sum_{t=1}^k\sum_{i=1}^m\,\beta_i^{t,m}\,q_i^{t,m}
\right\| & = &
\left\|\,\frac1k\sum_{t=1}^k\left(b-\sum_{i=1}^m\,\beta_i^{t,m}
\,q_i^{t,m}\right)\right\|
\\ \ \\
& = &
\left\|\,\frac1k\sum_{t=1}^k\sum_{i=1}^m\int_{Q_i^{t,m}}\,(b-\,\beta_i^{t,m})\,dE_{\cB}
\,\right\|
\\ \ \\
&\le&\left\|\,\frac1k\sum_{t=1}^k\sum_{i=1}^{n_1}\,\int_{Q_i^{t,m}}
\,(f-\beta_i^{t,m})\,dE_{\cB} \right\|+\epsilon
\\ \ \\
&\le&\left\|\,
\frac{2\,\|b\|}{k}\sum_{t=1}^k\sum_{i=1}^{n_1}\,q_i^{t,m}
\right\|+ \epsilon= \frac{2\|b\|}{k}+\epsilon\leq 2\epsilon
\end{eqnarray*}
where the first inequality is a consequence of 2' and the last
equality follows from 3'. \qed

\bigskip

\noindent {\it Proof of Lemma \ref{Adixmier}.} Fix a norm dense
subset $B=(b_j)_{j\in\NN}\subseteq \cB$. In the construction
leading to Dixmier's Theorem, a previous result
\cite[8.3.4]{kadison} asserts that for each $j$, there exists a
sequence $\{\rho_j^n\}_{n\in \NN} \subseteq\cD(\cM)$ such that for
every $1\leq h\leq j$,
$\|\rho_j^n(b_h)-\tau(b_h)\,I\|\xrightarrow[]{n}0$. For each $j\in
\NN$, let $n_0=n_0(j)\in \NN$ be such that if $n\geq n_0$ then
$\|\rho_j^n(b_h)-\tau(b_h)\,I\|\leq 1/j$ for $1\leq h\leq j$. If
we let $\rho_j=\rho_j^{n_0(j)}$ for $j\in \NN$,  we get
$\|\rho_j(b_h)-\tau(b_h)\,I\|\xrightarrow[]{j}0$ for every
 $h\in\NN$. Since $(b_j)_{j\in\NN}$ is norm dense in $\cB$
 we have $\lim_j\|\rho_j(b)-\tau(b)\,I\|=0$ for every $b\in\cB$.

For every $i=1,\ldots,m$, consider the factor $p_i\cM p_i$ with
(normalized) trace $\tau_i(p_ix)=\tau(xp_i)/\tau(p_i)\,$. By the
Dixmier approximation property mentioned in the first paragraph,
applied to the separable $C^*$-subalgebra $p_i\cB$ of the finite
factor $p_i\cM p_i$, there exists a sequence $\{\rho_j^i\}_{j\in
\NN}\in \cD(p_i\M p_i)$ such that $\lim_{j\rightarrow
\infty}\|\rho_j^i(p_ib)-\tau_i(p_ib)p_i\|=0$, for every $b\in\cB$.

For each $\rho\in \cD(p_i\cM p_i)$, we can consider an extension
$\tilde\rho\in\cD(\cM)$ as follows: if $\rho(p_i
b)=\sum_{h=1}^k\lambda_h\ u_h\,b\,u_h^*$, with $u_h\in \U(p_i\M
p_i)$, define $\tilde\rho\in \cD(\M)$ by
$\tilde\rho(b)=\sum_{h=1}^k\lambda_h\ \tilde u_h\,b\,\tilde
u_h^*$, where $\tilde u_h=u_h+(1-p_i)\in \U(\M)$. If $1\leq i\leq
m$ set $\rho_j=\prod_{i=1}^m \tilde \rho_j^{\,i}$ for $j\geq
1$. It is easy to verify that if $1\leq i\leq m$ then $
\rho_j(b\,p_i)= \tilde \rho_j^{\,i}(b\,p_i)$ for every $b\in\cB$.
Then, if $b\in\cB$,
\[\begin{array}{rcl}
\ds\left\|\rho_j(b)-\sum_{i=1}^m\beta_i(b)p_i\right\|
&=&
\ds\left\|\sum_{i=1}^m \tilde
\rho_j^{\,i}(b\,p_i)-\tau_i(b\,p_i)p_i\right\|
\xrightarrow[j\rightarrow \infty]{}0.\qed\end{array}
\]

\bigskip

\noindent{\bf Ackowledgements. }
We wish to thank Professors D. Farenick and D. Stojanoff for their support and useful discussions regarding the
material in this paper.

\end{document}